\documentclass[conference]{IEEEtran}
\IEEEoverridecommandlockouts
\usepackage{lcsys}
\usepackage{cite}
\usepackage{amsmath,amssymb,amsfonts}
\usepackage{algorithmic}
\usepackage{graphicx}
\usepackage{textcomp}

\usepackage{algorithmic}
\usepackage{graphicx}
\usepackage{textcomp}
\usepackage{caption, subcaption}
\usepackage{siunitx}
\usepackage{float}
\usepackage{gensymb}  
\usepackage{pstool,cleveref}
\usepackage{makecell}
\usepackage{yfonts}
\usepackage{multirow}
\usepackage{makecell}
\usepackage[inkscapelatex=false]{svg}
\usepackage{lipsum} 
\usepackage{geometry}

\newtheorem{remark}{Remark}
\newtheorem{prob}{Problem}
\newtheorem{property}{Property}

\def\BibTeX{{\rm B\kern-.05em{\sc i\kern-.025em b}\kern-.08em
    T\kern-.1667em\lower.7ex\hbox{E}\kern-.125emX}}
\begin{document}
\newgeometry{top=21.2mm, bottom=15.2mm, left=16.9mm, right=16.9mm}
\title{Scalable Optimal Motion Planning for Multi-Agent Systems by Cosserat Theory of Rods}
\author{Amirreza Fahim Golestaneh$^{1}$, \IEEEmembership{Member, IEEE}, Maxwell Hammond$^{1}$, and Venanzio Cichella$^{1}$, \IEEEmembership{Member, IEEE}
\thanks{This work was supported by NASA and ONR. }
\thanks{$^{1}$ The authors are with the Department of Mechanical Engineering, University of Iowa, Iowa City, IA 52240 {\tt\small \{amirreza-fahimgolestaneh, maxwell-hammond, venanzio-cichella\} @uiowa.edu} }
}

\maketitle

\begin{abstract}
We address the motion planning problem for large multi-agent systems, utilizing Cosserat rod theory to model the dynamic behavior of vehicle formations. The problem is formulated as an optimal control problem over partial differential equations (PDEs) that describe the system as a continuum. This approach ensures scalability with respect to the number of vehicles, as the problem's complexity remains unaffected by the size of the formation. The numerical discretization of the governing equations and problem's constraints is achieved through Bernstein surface polynomials, facilitating the conversion of the optimal control problem into a nonlinear programming (NLP) problem. This NLP problem is subsequently solved using off-the-shelf optimization software. We present several properties and algorithms related to Bernstein surface polynomials to support the selection of this methodology. Numerical demonstrations underscore the efficacy of this mathematical framework.
\end{abstract} 
\begin{IEEEkeywords}
Multi-agent system, Continuum mechanics, Cosserat theory of rods, Bernstein surface polynomial.
\end{IEEEkeywords}
%
\section{Introduction}
\label{Sec_Introduction}
Multi-agent systems surpass the capabilities of single-agent systems, providing essential advantages in various fields including search and rescue, surveillance, and exploration of hard-to-reach areas. 
Their impact extends across numerous sectors, revolutionizing practices in logistics and military operations. 
The significant promise of these systems has sparked extensive research into their control \cite{Xuan2002a,Choi2009a,Xiao2009b,Wang2012a,Qi2015a,Freudenthaler2016a}. 
Nevertheless, current control strategies for multi-agent systems still have potential for improvement, especially in addressing scalability issues. 

Scalability stands as a critical challenge in multi-agent systems' modeling and control, as issues like coordination and communication escalate with the increase of agent numbers. 
Given that the quantity of agents can vary dramatically based on the task and environment, methods which are independent of the number of necessary agents become extremely favorable. In recent literature, the use of PDE-based formulations in trajectory generation and tracking control has facilitated this independence by representing multi-agents as continuums deforming in time. 
In particular, the diffusion equation \cite{Qi2015a}, Burgers' equation \cite{Freudenthaler2016a,MEURER2011_a}, hyperbolic PDEs \cite{Selivanov2022_a}, parabolic distributed parameter systems \cite{FREUDENTHALER2020_a}, and principles from continuum mechanics \cite{MacLennan2019_a} are some examples of the methods used to represent multi-agents in this way. 
Intuitively, the PDEs chosen to represent a multi-agent system significantly impacts the behavior of the output from the path planner or controller. 
Ideally, this choice will be applicable to a wide range of agent dynamics while providing tools for safe constraint enforcement of agent interactions. To this end, the Cosserat theory of rods \cite{Antman_Book,Rubin_Book} is of interest given its flexibility. 
This model provides a mathematical framework for describing the deformation of slender elastic rods, and it is particularly noteworthy as it tracks position and orientation along space and time. 
These variables are readily transferable in meaning to multi-agent applications, describing position and heading of vehicles, and their derivatives translate to agent interactions along space, and vehicle dynamics in time. 

In \cite{hammond2023path} the authors implemented optimal-control based path planning algorithms for a soft robot modeled by Cosserat rod theory.
Bivariate Bernstein polynomials, also known as Bernstein surfaces, were utilized as a numerical tool to transform the problem into a nonlinear programming problem (NLP) dependent on the Bernstein coefficients, which can be solved with widely available optimization software \cite{kielas2019bebot,kielas2022bernstein}.
Within this context, the Bernstein basis provides numerically stable convergence and several desirable properties for safe constraint enforcement and obstacle avoidance, as demonstrated in \cite{kielas2022bernstein}, where the Bernstein basis is used for path planning of multiple vehicles, each modeled using ordinary differential equations.
In comparison to other polynomial bases used in pseudospectral methods, such as Legendre and Chebyshev, the Bernstein basis's unique properties ensure safety at any approximation order, enabling the potential for near real-time low order solutions without compromising constraint satisfaction \cite{cichella2017optimal,cichella2020optimal}.
Building on these ideas, we use Cosserat rods theory as the basis for our problem formulation and similarly employ Bernstein surfaces for solutions to the resulting problem. Cosserat rod theory, as opposed to traditional ODE modeling, allows us to treat the multi-agent system as a whole, overcoming scalability issues of conventional planning strategies.

Our primary contribution is to demonstrate the applicability of Cosserat theory in creating scalable motion planning algorithms for large vehicle formations. 
Section~\ref{Sec_ProblemFormulation} introduces kinematics PDEs for multi-agent motion using Cosserat theory, formulating motion planning and outlining constraints. 
Section~\ref{Sec_BernsteinApproximation} discusses discretizing continuum-based motion with Bernstein polynomials.
Section~\ref{Sec_ImplementationResultsAndDiscussions} presents numerical results, and Section~\ref{Sec_Conclusion} summarizes key findings and implications. 
%
\section{Problem Formulation}
\label{Sec_ProblemFormulation}
\subsection{Multi-agent system dynamics}
In standard Cosserat theory of rods \cite{Antman_Book}, the configuration of the rod is determined by the position vector $\mathbf{r}(s,t) \in \mathbb{R}^3$ of the centerline of the rod. 
The orientation of the rod cross section is described by the rotation matrix $R(s,t) \in \mathbb{SO}(3)$, where $s \in [0,s_f]$ denotes the curve length of the centerline of the undeformed (reference) rod and $t \in [0, t_f]$ represents time. The system of kinematic PDEs which defines the deformation of a Cosserat rod is given as
\begin{equation}
\begin{aligned}
    \mathbf{r}_s(s,t) &= R(s,t) \mathbf{l}(s,t), \\
    R_s(s,t) &= R(s,t) \mathbf{h}^\wedge(s,t), \\
\end{aligned}
\quad
\begin{aligned}
    \mathbf{r}_t(s,t)&= R(s,t) \mathbf{v}(s,t) 
    \\
    R_t(s,t)&= R(s,t) \boldsymbol{\omega}^\wedge(s,t) 
\end{aligned}
\label{Eq_KinematicsEqsOfStandardCosseratRod}
\end{equation}
where the operator $(\cdot)^\wedge$ returns the skew-symmetric matrix corresponding to a vector and subscripts $s$ and $t$ respectively denote the spatial and temporal derivatives.  
In these equations, the vectors $\mathbf{l}(s,t) = [\nu_1 , \nu_2 , \nu_3 ]^\top$ and $\mathbf{v}(s,t) = [ v_1 , v_2 , v_3 ]^\top$ represent the translational strain and velocity vectors, respectively. Similarly, $\mathbf{h}(s,t) = [ \mu_1 , \mu_2 , \mu_3 ]^\top$ and $\boldsymbol{\omega}(s,t) = [ \omega_1 , \omega_2 , \omega_3 ]^\top$ are vectors denoting the bending strain and angular velocity, respectively. 
\par
\begin{figure}[t]
\centerline{\includegraphics[width = 0.85\linewidth]{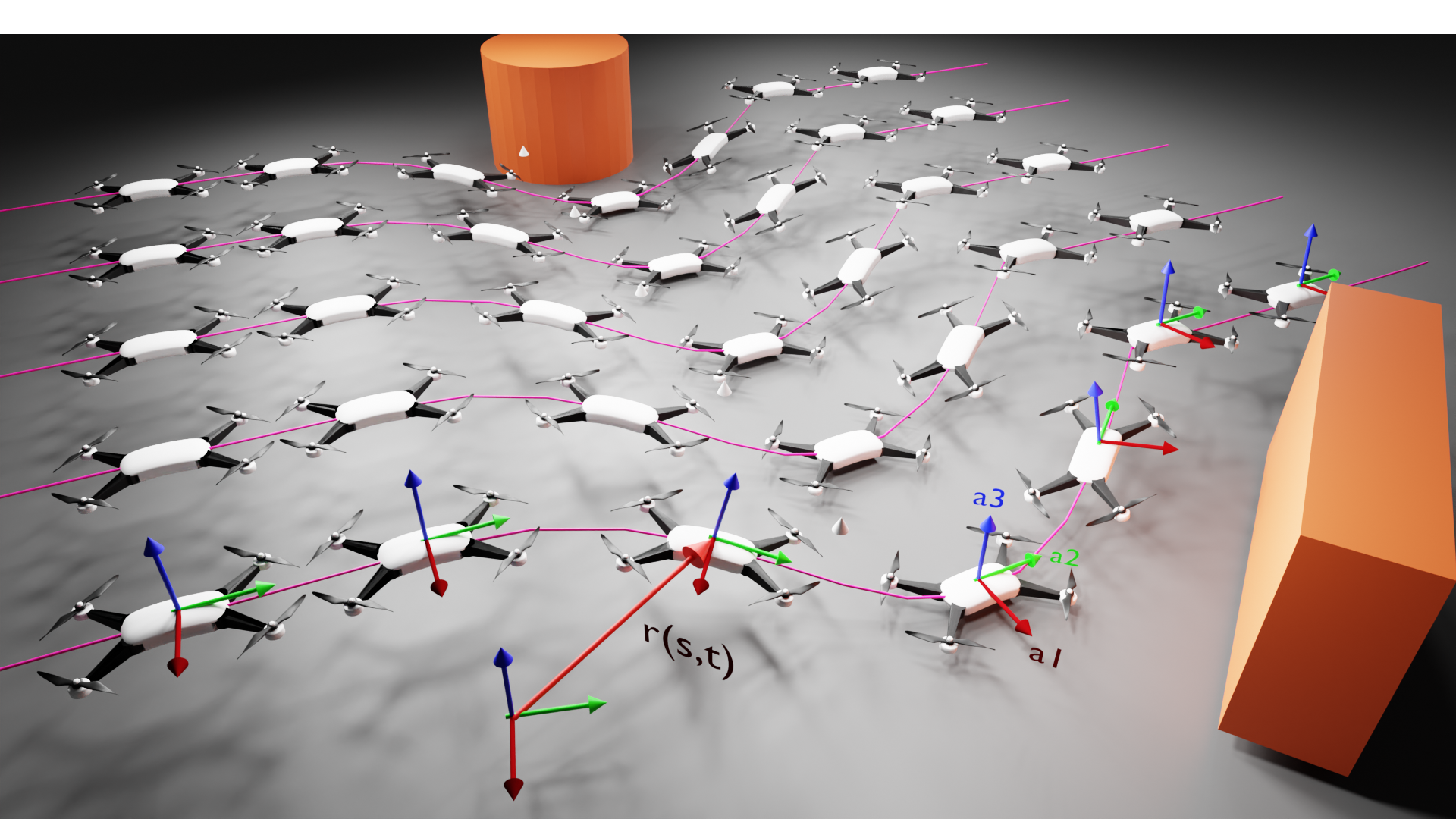} }
\caption{Multi-agent motion.} 
\label{Fig_Multi_Drones}
\end{figure}

Central to this paper is the idea that the variables of the Cosserat rod kinematics directly correspond to the pose of individual agents within the multi-agent system. The desired position for an agent at time $t$ corresponds to a specific point along the rod, denoted as $\mathbf{r}(s_i,t)$, where for equal spacing between agents $s_i=i \cdot s_f/n_v$, $i=1,\ldots,n_v$, and $n_v$ represents the number of vehicles. Similarly, $R(s_i,t)$ returns the orientation of a given agent in the system, and the temporal partial derivatives $\mathbf{v}(s_i,t)$ and $\boldsymbol{\omega}(s_i,t)$ give the vehicles linear and angular velocity.  Fig.\ref{Fig_Multi_Drones} provides a visual representation of this concept. Notice, by choosing $s_i$ as above, the problem's dependence on the number of vehicles  $n_v$ is removed, as increasing this parameter only changes how an individual agent's path is extracted from the solution. Consequently, only the rod's pose must be considered in problem formulation. The dynamic equations on the left-hand side of the kinematic system \eqref{Eq_KinematicsEqsOfStandardCosseratRod}, which represent spatial derivatives, are used to regulate interactions among vehicles within a fleet. Further, the equations on the right-hand side represent the dynamic constraints pertinent to each individual vehicle.

A noteworthy aspect of employing Cosserat rod theory in this context is the decoupling of temporal and spatial dynamics. For instance, in scenarios involving vehicles with more complex dynamic behavior, the standard system of equations can be substituted with a more general form:
\begin{equation*}
\begin{aligned}
    \mathbf{r}_s &= R \mathbf{l}, \\
    R_s &= R \mathbf{h}^\wedge, \\
\end{aligned}
\qquad
\begin{aligned}
    \mathbf{x}_t & = f(\mathbf{x},\mathbf{u})
    \\
    \mathbf{r} & = g_r(\mathbf{x},\mathbf{u}) \\
    \mathbf{R} & = g_R(\mathbf{x},\mathbf{u}) 
\end{aligned}
\end{equation*}
where $ {f}$, $ g_r$ and $ g_R$ represent the nonlinear motion dynamics specific to a given system, with \(\mathbf{x}\) and \(\mathbf{u}\) being some state vector and control input. Likewise, for more complex multi-agent configurations, the spatial derivative equations can also be adapted to embody more generalized interactions. However, the focus of this paper will remain on the system as described by \eqref{Eq_KinematicsEqsOfStandardCosseratRod}, thereby illustrating the application of Cosserat rod theory within the scope of multi-agent motion planning.

\subsection{Optimal motion planning}
\label{Sec_OptimalMotionPlanning}
Here, the motion planning problem for the multi-agent system is formulated as an optimal control problem over the PDE system in \eqref{Eq_KinematicsEqsOfStandardCosseratRod}. Subsequently, a detailed mathematical framework is introduced to outline the objectives that the motion planner is designed to fulfill.
%
\subsubsection{Feasibility Constraints}
\label{Sec_FeasibilityConstraintsForPoseBased}
Within the framework provided by \eqref{Eq_KinematicsEqsOfStandardCosseratRod}, the motion planning algorithm is tasked with generating trajectories for each vehicle that are dynamically feasible. The designed trajectories must enforce specific inter-vehicle constraints to ensure desired spacing. These constraints include preventing the vehicles from approaching too closely to one another or straying too far apart. Such requirements are enforced through the application of constraints on the variables of the Cosserat rod model. In particular, we consider the following constraints:
\begin{equation}
\begin{aligned}
\nu^2_{min} \leq & \Vert \mathbf{l}(s,t) \Vert^2 \leq \nu^2_{max} \\
& \Vert \mathbf{h}(s,t) \Vert^2  \leq \mu^2_{max} 
\end{aligned}
\qquad
\begin{aligned}
& \Vert \mathbf{v}(s,t) \Vert^2  \leq v^2_{max}  \\
& \Vert \boldsymbol{\omega}(s,t) \Vert^2  \leq \omega^2_{max}  
\end{aligned}
\label{Eq_StrainVelocityConstraintsPoseBased}
\end{equation}
for all \(s\in[0,s_f]\), \(t\in[0,t_f]\), where the norm \(||\cdot||\) denotes the Euclidean norm. The constraints on the left-hand side of \eqref{Eq_StrainVelocityConstraintsPoseBased} primarily address the inter-vehicle requirements by regulating both the position and orientation (attitude) of each vehicle with respect to its neighbour. 
For the sake of illustration, consider again the case of equally spaced agents located at $s_i$. The separation between neighboring agents at time \( t \) can be constrained by evaluating the ratio \( \frac{\mathbf{r}(s_{i+1},t) - \mathbf{r}(s_i,t)}{s_{i+1} - s_i} \). This ratio effectively represents the strain. Transitioning to the continuum case, this strain corresponds to \( \mathbf{l}(s,t) \). Analogous reasoning is applicable to the analysis of attitude dynamics.
Conversely, the constraints on the right-hand side are focused on the individual dynamics of each vehicle, ensuring that their speed and angular velocity profiles remain bounded. \par
%
%
\subsubsection{Initial and Final Conditions}
\label{Sec_InitialFinalConditionsPoseBased}
The motion planning algorithm must incorporate initial constraints to ensure the trajectory aligns with the agents' initial positions and attitudes. Final constraints can be included if the terminal pose needs to match a predefined spatial formation. These conditions can be formulated on the Cosserat rod states as follows:
\begin{align}
IC(\mathbf{r}(s,0),R(s,0), \mathbf{l}(&s, 0),\mathbf{h}(s,0), \mathbf{v}(s,0), \boldsymbol{\omega}(s,0))= 0 \nonumber \\[5pt] 
FC(\mathbf{r}(s,t_f),R(s,t_f), &\mathbf{l}(s,t_f),\mathbf{h}(s,t_f), \nonumber \\
& \mathbf{v}(s,t_f), \boldsymbol{\omega}(s,t_f))= 0.
\label{Eq_ICsFCsPoseBasedStrainBased}
\end{align}
%
\subsubsection{Obstacle Avoidance}
\label{Sec_ObstacleAvoidancePoseBased}
Finally, obstacle avoidance constraints need to be formulated to ensure that the optimal planned motion of the multi-agent system preserves a requisite safe distance from any obstacles within the operational workspace. Specifically, any point $\mathbf{r}(s,t)$ needs to remain at a distance greater than a specified safety margin, $\epsilon$, from any point $\mathbf{p}$ on the boundary $\partial P$ of an obstacle $P \subset \mathbb{R}^3$. This is formally represented as:
\begin{align}
d(\mathbf{r}(s,t), \mathbf{p}) > \epsilon, \quad \forall s\in[0,s_f], \, \forall t\in[0,t_f], \, \forall \mathbf{p} \in \partial P
\label{Eq_ObstacleAvoidanceConstraintFormula}
\end{align}
where $d(\mathbf{r}(s,t), \mathbf{p})$ denotes the distance between $\mathbf{r}(s,t)$ and the obstacle point $\mathbf{p}$. 

Finally, the motion planning problem for multi-agent system can be formulated as follows:

\begin{prob} 
\label{prob:multivehicle_continuous}
Find \(s_f, t_f, \mathbf{r},R, \mathbf{l},\mathbf{h}, \mathbf{v},\) and \(\boldsymbol{\omega}\) that minimize
\begin{equation}
    \int_0^{s_f} \int_0^{t_f} 
    {\ell}(\mathbf{r},R,\mathbf{l},\mathbf{h},\mathbf{v}, \boldsymbol{\omega}) dt ds.
    \label{Eq_GeneralCostFunction}
\end{equation}
subject to \eqref{Eq_KinematicsEqsOfStandardCosseratRod}, \eqref{Eq_StrainVelocityConstraintsPoseBased}, \eqref{Eq_ICsFCsPoseBasedStrainBased}, and \eqref{Eq_ObstacleAvoidanceConstraintFormula}.
\end{prob}

Here, the running cost function \(\ell()\) can be tailored to reflect goals pertinent to the collective operation of the multi-agent system. These objectives may vary, ranging from spatial considerations—such as maintaining proximity—to temporal goals, including minimizing travel time. If needed, the cost can be modified to incorporate terminal costs. However, in most cases of interest, the proposed formulation is general enough.

%
\section{Bernstein Surface Approximation }
\label{Sec_BernsteinApproximation}
In this section, we explore the application of bivariate Bernstein polynomials as a numerical approximation technique for addressing the motion planning problem at hand. By representing the variables of interest as Bernstein surfaces, Problem \ref{prob:multivehicle_continuous} is reformulated into a NLP, which can be solved using off-the-shelf optimization software. This approach was initially introduced in \cite{hammond2023path} for the motion planning of soft robots. In this paper, we extend the concept to develop a scalable algorithm for multi-agent systems planning.  

\begin{figure}
   \centering
   \includegraphics[width = 1\linewidth,trim={0.5cm 0.3cm 0.9cm 0.4cm},clip]{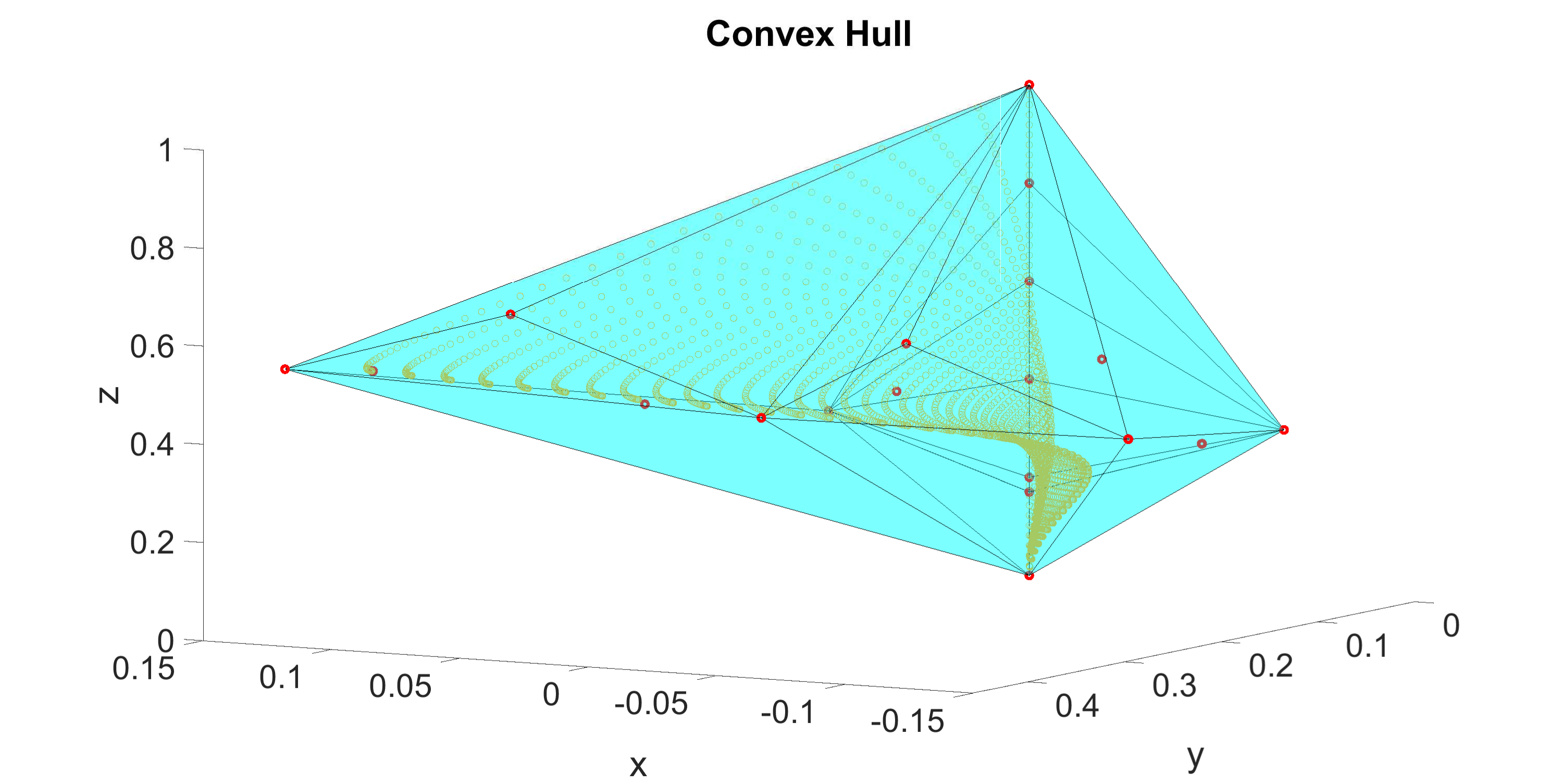}
   \caption{Bernstein surface of order \(m=5, n=5\). The red dots represent control points, the orange dots represent points on the surface, and the blue polygon is the convex hull.}
   \label{fig:convex_hull}
\end{figure}

Recall that the desired position assigned to agent \(i\) at time $t$ is denoted as $\mathbf{r}(s_i,t)$.
We let $\mathbf{r}(s,t)$, \(s\in[0,s_f]\), \(t\in[0,t_f]\) be parameterized by an \(m \times n \) Bernstein surface defined as 
\begin{equation} \label{eq:BSposition}
    \quad \textbf{r}(s,t) = \sum_{i=0}^{m}\sum_{j=0}^{n}\bar{\textbf{r}}_{i,j}^{m,n}B_{i}^{m}(s)B_{j}^{n}(t).
\end{equation}
where $\bar{\textbf{r}}_{i,j}^{m,n}\in \mathbb{R}^{3}$, \(i=0,\ldots,m\), \(j=0,\ldots,n\) are control points and $B_{i}^{m}(s)$ is the Bernstein polynomial basis over the interval $[0,s_f]$,
\begin{equation*} 
    B_{i}^{m}(s)={m \choose i}\frac{s^i\left(s_f-s\right)^{m-i}}{s_f^m}, \quad 0 \le i \le m, 
\end{equation*}
where ${m \choose i} = \frac{m!}{i!\left(m-i\right)!}.$ Fig.~\ref{fig:convex_hull} depicts an example of Bernstein surface (orange) defined by its control points (red circles).
Similarly, consider the desired orientation assigned to the \(i\)th agent at time \(t\), which can be characterized by the Euler angles \((\phi(s_i,t), \theta(s_i,t), \psi(s_i,t))\). For the angle \(\phi(s,t)\), the Bernstein surface is defined by:
\begin{equation} \label{eq:phi}
    \phi(s,t) = \sum_{i=0}^{m}\sum_{j=0}^{n}\bar{\phi}_{i,j}^{m,n} B_{i}^{m}(s) B_{j}^{n}(t),
\end{equation}
where \(\bar{\phi}_{i,j}^{m,n}\) denotes the control points associated with the surface. The angles \(\theta(s,t)\) and \(\psi(s,t)\) are analogously defined, utilizing their respective control points \(\bar{\theta}_{i,j}^{m,n}\) and \(\bar{\psi}_{i,j}^{m,n}\), for \(i=0,\ldots,m\) and \(j=0,\ldots,n\). Similarly, the variables \(\mathbf{l}(s,t)\), \(\mathbf{h}(s,t)\), \(\mathbf{v}(s,t)\), and \(\boldsymbol{\omega}(s,t)\) are characterized by equivalent summations, employing their respective control points, \(\bar{\mathbf{l}}_{i,j}^{m,n}\), \(\bar{\mathbf{h}}_{i,j}^{m,n}\), \(\bar{\mathbf{v}}_{i,j}^{m,n}\), and \(\bar{\boldsymbol{\omega}}_{i,j}^{m,n}\) for \(i=0,\ldots,m\) and \(j=0,\ldots,n\).

This parameterization allows the transcription of the motion planning problem for multivehicle systems into an optimization problem utilizing the control points of the aforementioned variables as the optimization variables:
\begin{prob} 
\label{prob:multivehicle_discrete}
Find \(s_f\), \(t_f\), \(\bar{\textbf{r}}_{i,j}^{m,n}\), \(\bar{\phi}_{i,j}^{m,n}\), \(\bar{\theta}_{i,j}^{m,n}\), \(\bar{\psi}_{i,j}^{m,n}\), \(\bar{\textbf{l}}_{i,j}^{m,n}\), \(\bar{\textbf{h}}_{i,j}^{m,n}\), \(\bar{\textbf{v}}_{i,j}^{m,n}\) and \(\bar{\boldsymbol{\omega}}_{i,j}^{m,n}\), \(i=0,\ldots,m\) and \(j=0,\ldots,n\), that minimize
\begin{equation}
    \sum_{i=0}^{m} w_{s,i} \sum_{j=0}^{n} w_{t,j} 
    \tilde{\ell}(\bar{\textbf{r}}_{i,j}^{m,n},R_{ij}, \bar{\textbf{l}}_{i,j}^{m,n}, \bar{\textbf{h}}_{i,j}^{m,n}, \bar{\textbf{v}}_{i,j}^{m,n}, \bar{\boldsymbol{\omega}}_{i,j}^{m,n}) 
    \label{Eq_GeneralCostFunction_discrete}
\end{equation}
subject to \eqref{Eq_KinematicsEqsOfStandardCosseratRod}, \eqref{Eq_StrainVelocityConstraintsPoseBased}, \eqref{Eq_ICsFCsPoseBasedStrainBased}, and \eqref{Eq_ObstacleAvoidanceConstraintFormula}.
\end{prob}

In \eqref{Eq_GeneralCostFunction_discrete}, \(R_{ij}\) is defined as \(R_{ij} = R(\bar{\phi}_{i,j}^{m,n},\bar{\theta}_{i,j}^{m,n}, \bar{\psi}_{i,j}^{m,n}) . \)
The cost function is approximated using numerical integration over the Bernstein polynomial basis, with the weights \(w_{s,i} = \frac{s_f}{m+1}\) and \(w_{t,i} = \frac{t_f}{n+1}\) for all \(i=0,\ldots,m\), \(j=0,\ldots,n\). Next, we present algorithms and properties of Bernstein surfaces that facilitate the approximation of the constraints. These underscore the main advantages of employing this particular approximation method over alternative approaches.

\begin{property}[Arithmetic operations]
Addition and subtraction between two Bernstein surfaces can be performed directly through the addition and subtraction of their control points. The control points of the Bernstein surface \(y(\cdot,\cdot)\) resulting from multiplication between two Bernstein surfaces, \(g(\cdot,\cdot)\) and \(h(\cdot,\cdot)\) with control points $\bar{g}_{i,j}^{m,n}$ and $\bar{h}_{k,l}^{a,b}$ can be obtained by 
\begin{equation}
    \begin{split}
        \bar{y}_{e,f}^{m+a,n+b} &= \\ \sum_{q=max(0,e-a)}^{min(m,e)}&\sum_{r=max(0,f-b)}^{min(n,f)}  
        \frac{{m \choose q}{n \choose r}{a \choose e-q}{b \choose f-r}}{{m+a \choose e}{n+b \choose f}}\bar{g}_{q,r}^{m,n}\bar{h}_{e-q,f-r}^{a,b}
    \end{split}
\end{equation}
\end{property}

\begin{property}[Derivatives]
The partial derivatives of a Bernstein surface can be calculated by multiplying a differentiation matrix with the surface's control points. For example, consider the Bernstein surface representing the rotation of the rod, \eqref{eq:phi}, and let \(\bar{\mathbf{\Phi}}^{m,n}\) be the matrix of control points, i.e., \(\{\bar{\mathbf{\Phi}}^{m,n}\}_{i,j} = \bar{\phi}_{i,j}^{m,n}\). The partial derivatives $\frac{\partial}{\partial s} \phi(s,t)$ and $\frac{\partial}{\partial t} \phi(s,t)$ are given by control points
\begin{equation} \label{eq:bsurfderiv}
        \bar{\mathbf{\Phi}}^{m,n}_s=\textbf{D}_m^\top \bar{\mathbf{\Phi}}^{m,n}, \quad
        \bar{\mathbf{\Phi}}^{m,n}_t=\bar{\mathbf{\Phi}}^{m,n}\textbf{D}_n ,
\end{equation}
where $\textbf{D}_m$
and $\textbf{D}_n$ are square differentiation matrices, see \cite{kielas2022bernstein}. 
\end{property}

Using the two properties above, the dynamic constraints \eqref{Eq_KinematicsEqsOfStandardCosseratRod} can be expressed as a system of algebraic equations. 
The expressions in \eqref{Eq_StrainVelocityConstraintsPoseBased} can be represented by
\begin{equation} \label{eq:vexpressions}
\begin{aligned}
||\mathbf{l}(s,t)||^2 & = \sum_{i=0}^{2m}\sum_{j=0}^{2n}\bar{{l}}_{i,j}^{2m,2n}B_{i}^{2m}(s)B_{j}^{2n}(t) \\ 
||\mathbf{v}(s,t)||^2 & = 
\sum_{i=0}^{2m}\sum_{j=0}^{2n}\bar{{v}}_{i,j}^{2m,2n}B_{i}^{2m}(s)B_{j}^{2n}(t) \\
||\, \mathbf{h}(s,t) \, ||^2 & = 
\sum_{i=0}^{2m}\sum_{j=0}^{2n}\bar{{h}}_{i,j}^{2m,2n}B_{i}^{2m}(s)B_{j}^{2n}(t) \\
|| \, \boldsymbol{\omega}(s,t) \, ||^2 & = 
\sum_{i=0}^{2m}\sum_{j=0}^{2n}\bar{{\omega}}_{i,j}^{2m,2n}B_{i}^{2m}(s)B_{j}^{2n}(t) \\
\end{aligned}
\end{equation}
where the coefficients $\bar{{l}}_{i,j}^{2m,2n}$, $\bar{{v}}_{i,j}^{2m,2n}$, $\bar{{h}}_{i,j}^{2m,2n}$, and \(\bar{{\omega}}_{i,j}^{2m,2n}\) can be obtained from  algebraic manipulation of the Bernstein coefficients of $\mathbf{l}(s,t)$, $\mathbf{v}(s,t)$, $\mathbf{h}(s,t)$ and $\boldsymbol{\omega}(s,t)$.

\begin{figure}
   \centering
   \includegraphics[width = 1\linewidth,trim={0.5cm 0.3cm 0.9cm 0.4cm},clip]{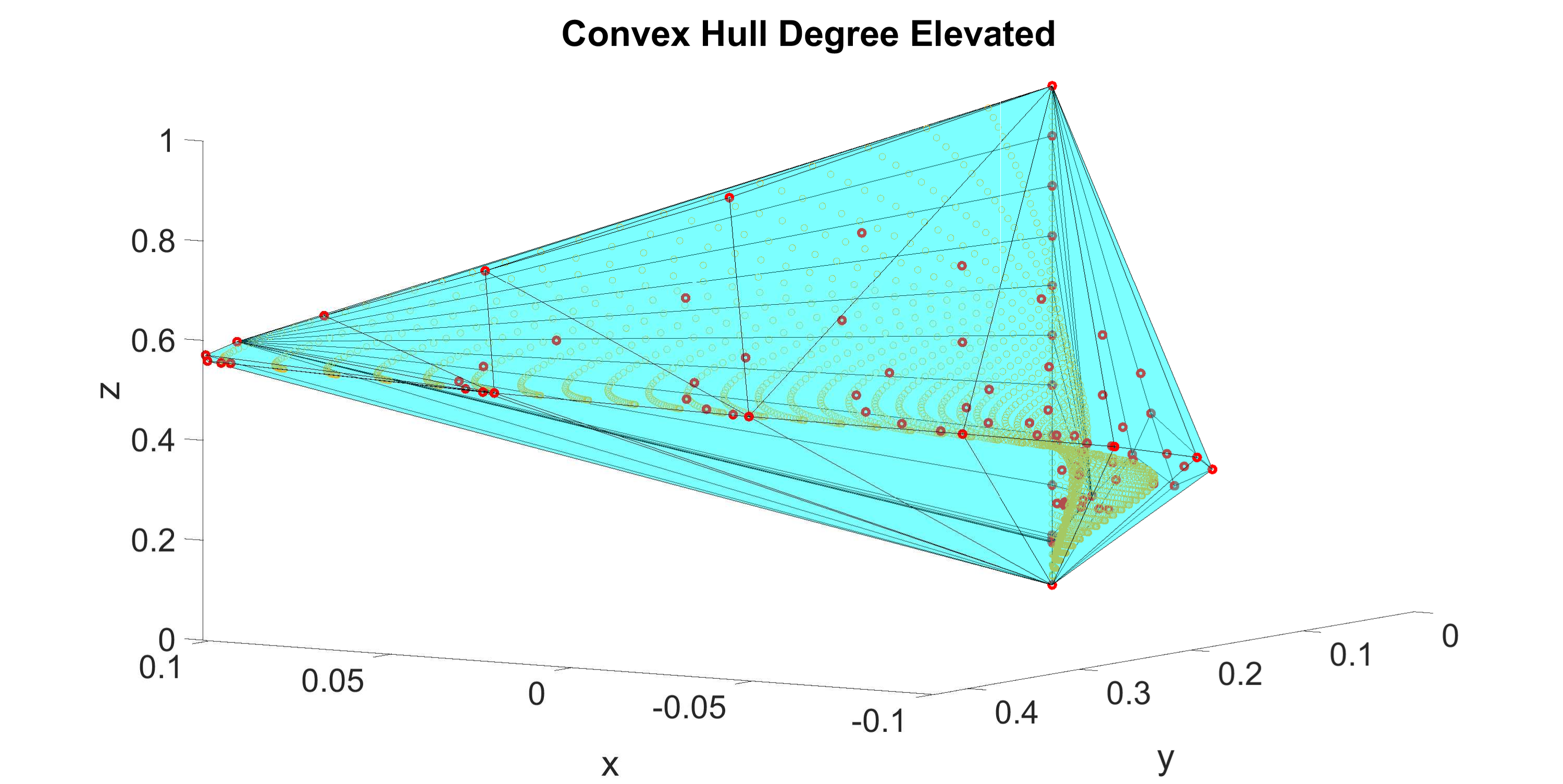}
   \caption{Bernstein surface from Fig.~\ref{fig:convex_hull} after being degree elevated to \(m=10, n=10\).}
   \label{fig:deg_el}
\end{figure}

\begin{figure}
   \centering
   \includegraphics[width = 1\linewidth,trim={0.5cm 0.2cm 0.9cm 0.4cm},clip]{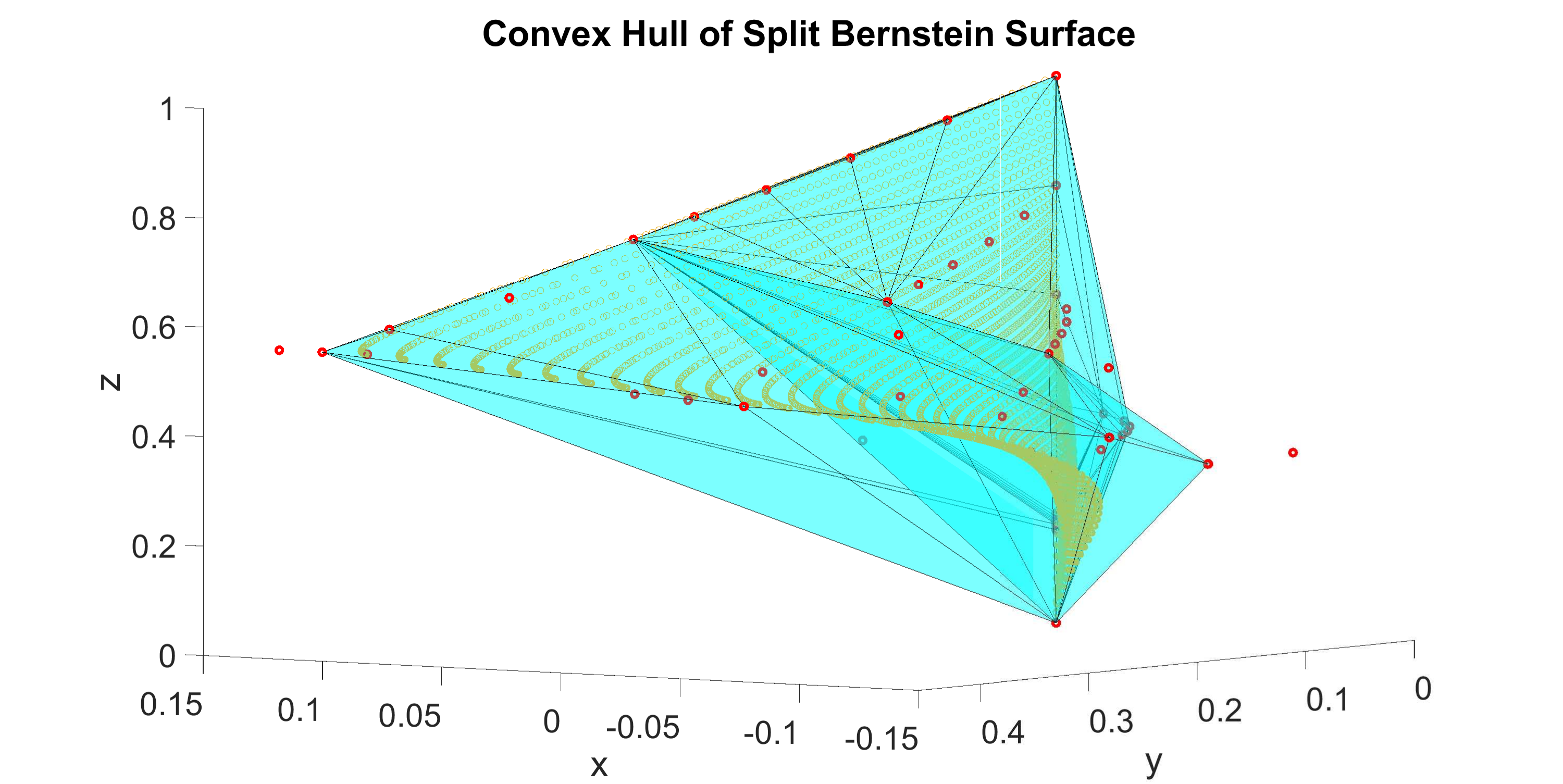}
   \caption{Bernstein surface from Fig.~\ref{fig:convex_hull} being split along the red line using the de Casteljau algorithm.}
   \label{fig:deCastel}
\end{figure}

To enforce the inequality constraints we use the convex hull property of Bernstein surfaces.
\begin{property}[Convex Hull]
A Bernstein surface lies within the convex hull defined by its control points. 
For example, the surface described in \eqref{eq:vexpressions} satisfies:
\begin{equation}
\min_{i,j} {\bar{{v}}_{i,j}^{2m,2n}}
\leq
||\mathbf{v}||^2 
\leq 
\max_{i,j} {\bar{{v}}_{i,j}^{2m,2n}}
\end{equation}
where $i=0,\ldots,2m$ and $j=0,\ldots,2n$.
\end{property}
This property can be used to impose boundary constraints by directly applying them to the control points of the surface. In fact, the dynamic limits in \eqref{Eq_StrainVelocityConstraintsPoseBased} can be imposed as follows
\begin{equation} \label{eq:vconstraintsdisc}
\begin{aligned}
&\nu_{\min}^2 \leq \bar{{l}}_{i,j}^{2m,2n} \leq \nu_{\max}^2, \quad &\bar{{v}}_{i,j}^{2m,2n} \leq v_{\max}^2, \\
&\bar{{h}}_{i,j}^{2m,2n}\leq \mu_{\max}^2, \quad &\bar{{\omega}}_{i,j}^{2m,2n}\leq \omega_{\max}^2, 
\end{aligned}
\end{equation}
for all \(i=0,\ldots,2m\), \(j=0,\ldots,2n\).

\begin{remark}
The convex hull that encloses the Bernstein surface may be significantly larger that the surface itself, leading to conservativeness when imposing bounds. To mitigate this conservativeness, degree elevation (see Fig.~\ref{fig:deg_el}) and the de Casteljau split (see Fig.~\ref{fig:deCastel}) can be employed. These are numerically stable methods that reduce the convex hull's dimensions, aligning it more closely with the actual surface. For a detailed description on these techniques and their application, the reader is directed to \cite{kielas2022bernstein}. 
\end{remark}

Regarding the boundary conditions, the endpoint value property of Bernstein surfaces is relevant.
\begin{property}[End point values]
The terminal points of a Bernstein surface coincide with their corresponding terminal control points. For instance,
\begin{equation}
\begin{split}
    \textbf{r}(0,0) = \bar{\textbf{r}}^{m,n}_{0,0}, \quad \textbf{r}(0,t_f) = \bar{\textbf{r}}^{m,n}_{0,n}, \quad 
    \textbf{r}(s_f,0) = \bar{\textbf{r}}^{m,n}_{m,0}, 
\end{split}
\end{equation}
Moreover, the surface's edges can be represented by the Bernstein polynomials of the edge's control points. For example,
\begin{equation}
    \textbf{r}(s,0) = \sum_{i=0}^m \bar{\textbf{r}}^{m,n}_{i,0} B_{i}^{m}(s).
\end{equation}
\end{property}
Leveraging this property, the boundary conditions specified in \eqref{Eq_ICsFCsPoseBasedStrainBased} can be expressed as functions of the control points of the Bernstein surfaces corresponding to the variables of interest.

To address obstacle avoidance constraints, the minimum distance algorithm for Bernstein surfaces is utilized \cite{kielas2022bernstein, chang2011computation}. This algorithm calculates the minimum distance between a Bernstein surface and a convex shape, corresponding to the left-hand side of \eqref{Eq_ObstacleAvoidanceConstraintFormula}. It integrates the Convex Hull property, the Endpoint Values property, the de Casteljau Algorithm, and the Gilbert-Johnson-Keerthi (GJK) algorithm \cite{gilbert1988fast}, combining these techniques to efficiently determine the minimum separation distance.
We direct the reader to \cite{kielas2022bernstein} for further discussion.
The inputs to the algorithm are the sets of Bernstein coefficients defining the Bernstein surface representing the multi-agent system moving in time, i.e., $\bar{\textbf{r}}^{m,n}$, and the vertices of a convex shape that must be avoided, i.e., $\bar{\textbf{Q}}$. The function provides the minimum distance between the two objects as its result, i.e., $\min_{s,t} d(\mathbf{r}(s,t),\mathbf{p})$, where \(d()\) was defined in \eqref{Eq_ObstacleAvoidanceConstraintFormula}.

By leveraging the properties of the Bernstein basis, Problem \ref{prob:multivehicle_discrete} can be solved as a NLP over Bernstein polynomial coefficients. When optimal control problems are approximated into NLP using Bernstein polynomials, the solutions to the NLP converge uniformly to the optimal solution of the original control problem as the order of approximation increases \cite{cichella2019consistent}. Of course, increasing the order leads to an increased computational cost. However, guarantees on constraint satisfaction afforded by the unique properties of the chosen basis offer a balance between computational efficiency and optimality without compromising on safety. By extension, the end user can determine an appropriate order for their application. For example, when an efficient solution is needed immediately to react to an unforeseen obstacle, a low-order safe solution can be generated almost instantly. Alternatively, a solution that closely approximates the optimal one can be computed in advance for a mission. The computational efficiency of this approach is furthered in the context of this work by the independence of this approach from the number of agents. 

%
%

%

    \begin{figure}[t]
    \centering
        \includegraphics[width=0.40\textwidth,trim={1cm -0.1cm 1cm 1cm},clip]{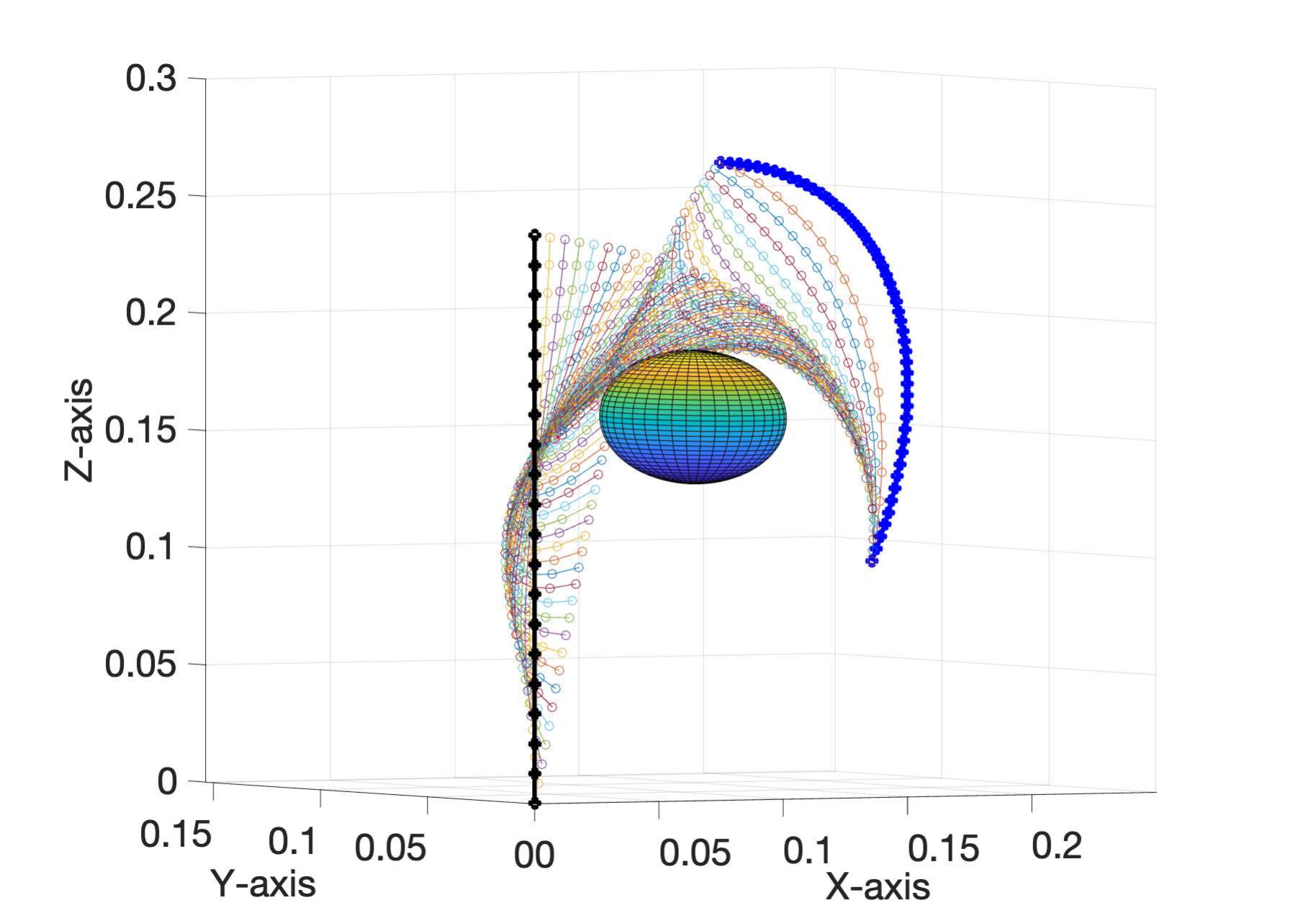}
        \caption{Case 1 - Optimal motion of agents from initial to final configurations. The black straight line depicts the initial formation of agents and the blue curve is a quarter of an ellipse for final formation, independent of number of agents. The colored curves depict the transitions of agents from their initial to final formations. }
        \label{Fig_7a_1_Ellipse_WithAvoidance}
    \end{figure}
    %
    %
    \begin{figure}[t]
    \centering
        \includegraphics[width = 0.70\linewidth,trim={1.3cm -0.1cm 1.2cm -0.1 cm},clip]{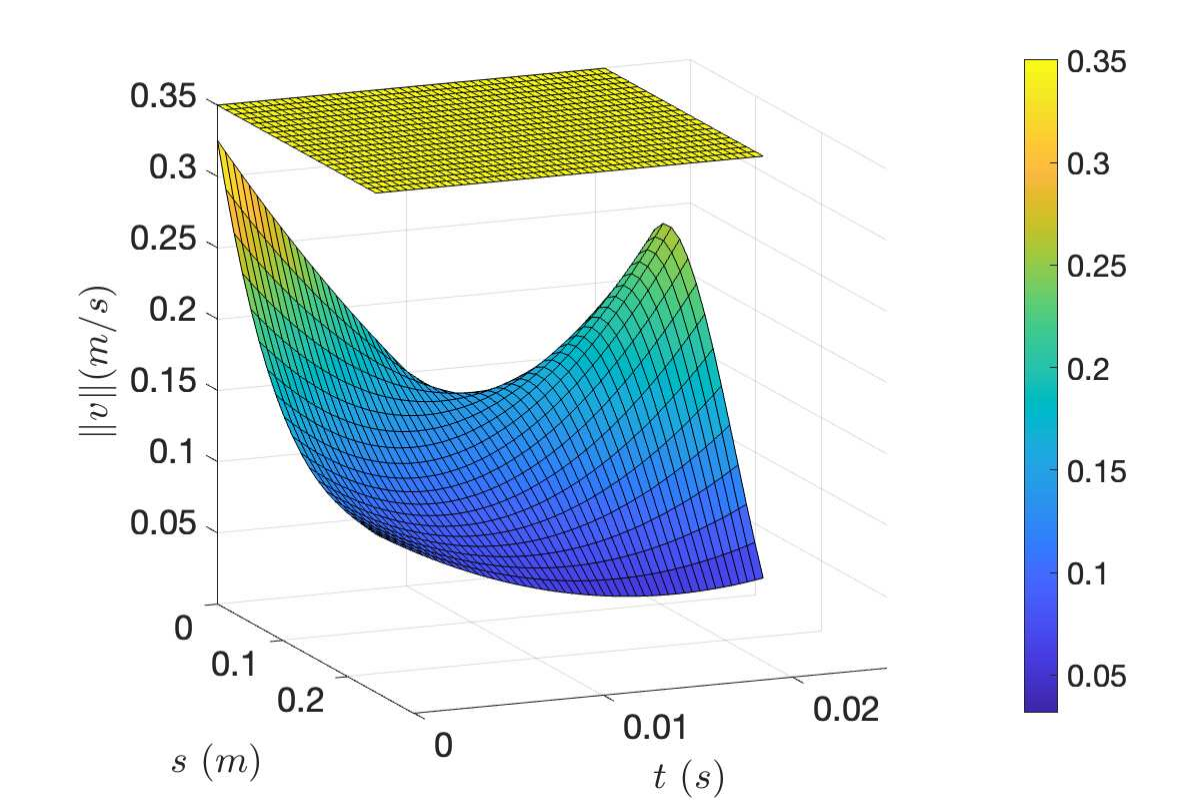}
        \caption{Case 1 -Translational velocities comply with constraints.}
        \label{Fig_8_1_Ellipse_WithAvoidance_VelocityVNorm}
\end{figure}
    %
%
\begin{figure}[t]
    \centering
        \includegraphics[width = 0.65\linewidth,trim={0.9cm 0.3cm 0.9cm 1cm},clip]{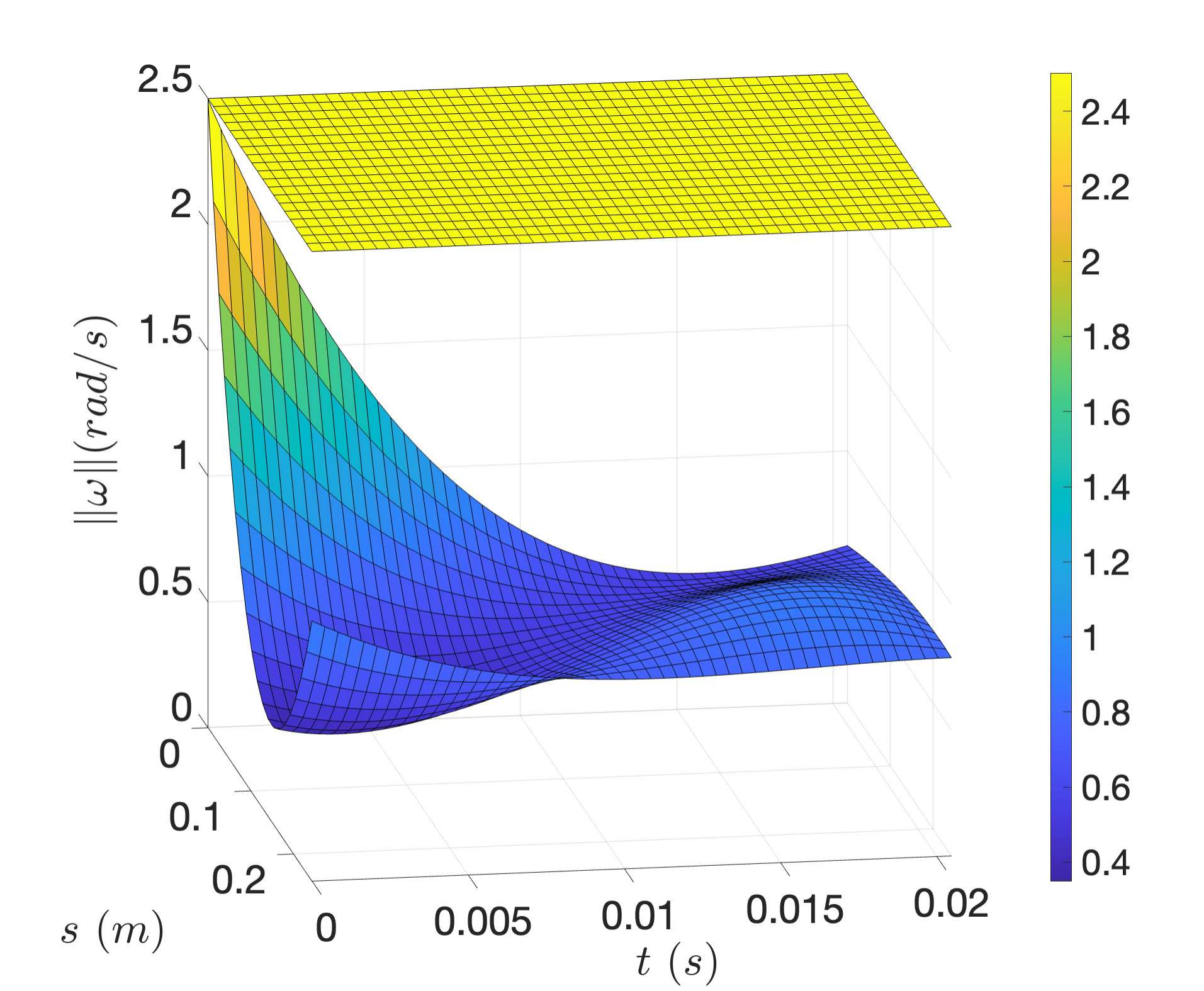}
        \caption{ Case 1 - Angular velocities comply with constraints.}
        \label{Fig_8_1_Ellipse_WithAvoidance_VelocityOmegaNorm}
\end{figure}
    \begin{figure}[t]
    \centering
        \includegraphics[width=0.55\textwidth,trim={0cm 0cm 0cm 0cm},clip]{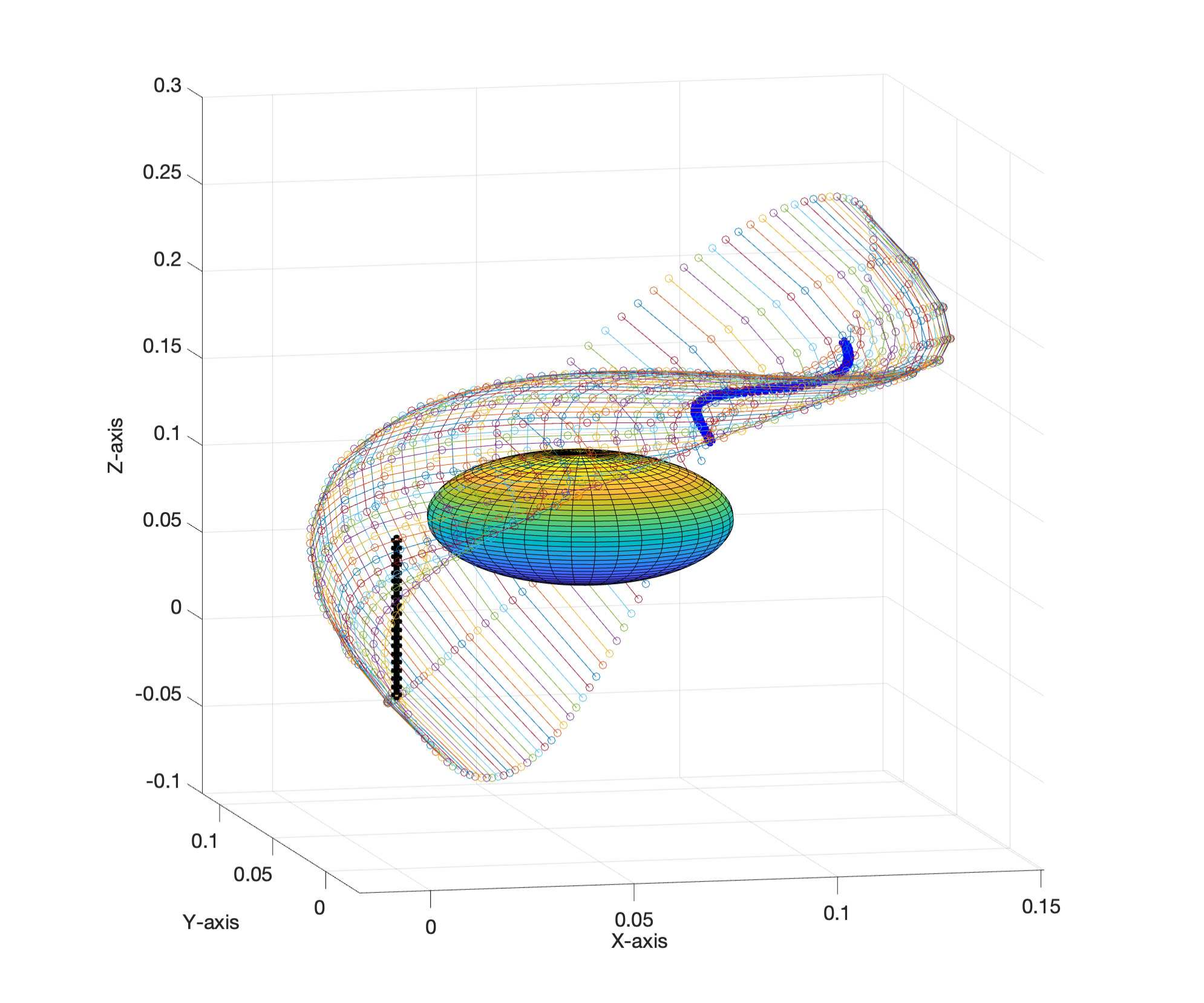}
        \caption{Case 2 - Optimal motion of agents from initial and final configuration. The black straight line depicts the initial formation of agents and the blue curve is a helix for final formation, independent of number of agents. The colored curves depict the transitions of agents from their initial to final formations.}
        \label{Fig_7a_1_Helix_WithAvoidance}
    \end{figure}
    %
%
%
\section{Numerical Results}
\label{Sec_ImplementationResultsAndDiscussions} 
%
In this section we discuss two cases of multi-agent motions with elliptical and helical final formations.
The cost function for these specified scenarios is defined as follows:

\begin{align}
 \int_{0}^{t_f} \{ & \Vert \mathbf{r}(0,t) - \mathbf{r}^\circledast_{des} \Vert_2^2 + \Vert \boldsymbol{\theta}(0,t) - \boldsymbol{\theta}^\circledast_{des} \Vert_2^2 + \nonumber\\[5pt] 
+ & \Vert \mathbf{r}(s_f,t) - \mathbf{r}^f_{des} \Vert_2^2 + \Vert \boldsymbol{\theta}(s_f,t) - \boldsymbol{\theta}^f_{des} \Vert_2^2 \} dt.
\label{Eq_CostFunctionPose}
\end{align}
In this formulation, \(\mathbf{r}(0,t)\) and \(\mathbf{r}(s_f,t)\) represent the positions of two agents at the formation's extremities, namely the leaders, with analogous definitions for the angles \(\boldsymbol{\theta}(0,t)\) and \(\boldsymbol{\theta}(s_f,t)\). 
The cost function's goal is to ensure that these two vehicles achieve convergence to the desired positions and orientations, denoted by \(\mathbf{r}^\circledast_{des}\), \(\boldsymbol{\theta}^\circledast_{des}\), \( \mathbf{r}^f_{des}\) and \(\boldsymbol{\theta}^f_{des}\), in the shortest possible time. 
The problem is approximated using a Bernstein surface of order \(m=6, n=6\), providing satisfactory
accuracy for our results. The optimization problem is executed in MATLAB, using fmincon optimization. The initial guess used to start the optimizer is a static point surface with constant attitude. Given the constraints of the problem, discussed below, this initial guess is unfeasible. This demonstrates the capability of the approach to generate feasible and near-optimal solutions even when the initial guess is far from being ideal.

Table \ref{Table_CaseStudiesOfMotionPlanning} presents two case studies, outlining data relevant to each enforced constraint. 
In Case 1, the multi-agent system transitions from a straight line formation to an elliptical one, with the final formation representing a quarter of an ellipse oriented in three-dimensional space. 
Details regarding the ellipse's center point, dimensions, and orientation are specified in Table~\ref{Table_CaseStudiesOfMotionPlanning}.
%
\begin{table}[tbp]
\caption{Case studies of motion planning.}
\begin{center}
\begin{tabular}{|c|c|}
\hline
 \textbf{Case 1}& \textbf{Case 2} \\
\hline
 \makecell[l]{ Initial Formation: \textbf{line} \\ $\mathbf{X}_o= s \mathbf{e}_3$ \\ $s \in [0,0.24] \, m$ }  & \makecell[l]{ Initial Formation: \textbf{ line} \\ $\mathbf{X}_o= s \mathbf{e}_3$ \\ $s \in [0,0.24] \, m$ } \\
\hline
  \makecell[l]{ Final Formation: \textbf{Ellipse} \\ Center: $[0.1,0.1,0.1]^T \, m$ \\Semi-axis: $0.2, 0.1 \, m$ \\ Axis 1: $[0.53,0.26,0.80]^T \, m$ \\ Axis 2: $[0.45,-0.89,0]^T \, m$ \\ Parameterization: $ \theta \in [0, \pi/2] $ }   &  \makecell[l]{ Final Formation: \textbf{Helix} \\ Radius: $0.01 \, m$ \\ Pitch: $0.02 \pi \, m$ \\ Slope: $\pi/4$  Rad   } \\
\hline
 \makecell[l]{ Obstacle: \textbf{Sphere} \\ Center: $[0.11,0.05,0.16]^T \, m$ \\ Radius: $0.03 \, m$ } &  \makecell[l]{ Obstacle: \textbf{Sphere:} \\ Center: $[0.06,0.05,0.08]^T \, m$ \\ Radius: $0.04 \, m$ } \\
\hline
 \makecell[l]{ $\nu_{max}= 2.25$, $\nu_{min}= 0.7$  \\ $\mu_{max}= 1.55 $  $Rad$ \\ $v_{max}= 0.35$ $m/s$ \\ $\omega_{max}= 2.5$ $ Rad/s $ }  &  \makecell[l]{ $\nu_{max}= 2.25$, $\nu_{min}= 0.7$  \\ $\mu_{max}= 1.55 $  $Rad$ \\ $v_{max}= 0.35$ $m/s$ \\ $\omega_{max}= 2.5$ $ Rad/s $ }  \\
\hline
\end{tabular}
\label{Table_CaseStudiesOfMotionPlanning}
\end{center}
\end{table}
%
%
The output of Case 1 is summarized in the Figs. \ref{Fig_7a_1_Ellipse_WithAvoidance}-\ref{Fig_8_1_Ellipse_WithAvoidance_VelocityOmegaNorm}. In Fig.  \ref{Fig_7a_1_Ellipse_WithAvoidance} black lines represent the initial and final formation of the agents, and the continuous manifold between these lines is the optimal trajectory of the system according to the scenario detailed in Table \ref{Table_CaseStudiesOfMotionPlanning}. Figs. \ref{Fig_8_1_Ellipse_WithAvoidance_VelocityVNorm} and \ref{Fig_8_1_Ellipse_WithAvoidance_VelocityOmegaNorm} represent the satisfaction of linear and angular velocity constraints, respectively. In each Fig., the upper bound of the relevant variable is shown as a plane over \(s\in[0,s_f]\) and \(t\in[0,t_f]\), and the plotted surface demonstrates adherence to this bound within the domain. Similar plots for constraints on spatial partial derivatives are omitted for brevity. The output of Case 2 is similarly represented in Fig. \ref{Fig_7a_1_Helix_WithAvoidance}.
\section{Conclusion}
\label{Sec_Conclusion}
This paper presented a scalable, PDE-based motion planning framework for large multi-agent systems, utilizing the Cosserat theory of rods. By employing Bernstein surface polynomials for discretization, we transformed the continuous problem into a solvable NLP, demonstrating the efficacy of our approach through numerical results. Future research will explore alternative continuum mechanics models tailored to specific operational formations. For instance, while the Cosserat rod model excels in planning for vehicles flying in 1D formation, two-dimensional formations may benefit from an extended Cosserat rod model or a Cosserat shell approach.

\bibliographystyle{plain} 
\bibliography{MultiAgent}

\begin{thebibliography}{10}

\bibitem{Antman_Book}
Stuart~S. Antman.
\newblock {\em Nonlinear Problems of Elasticity}.
\newblock Springer, 2005.

\bibitem{chang2011computation}
Jung-Woo Chang, Yi-King Choi, Myung-Soo Kim, and Wenping Wang.
\newblock Computation of the minimum distance between two b{\'e}zier curves/surfaces.
\newblock {\em Computers \& Graphics}, 35(3):677--684, 2011.

\bibitem{Choi2009a}
Jongeun Choi, Songhwai Oh, and Roberto Horowitz.
\newblock Distributed learning and cooperative control for multi-agent systems.
\newblock {\em Automatica}, 45(12):2802--2814, 2009.

\bibitem{cichella2017optimal}
Venanzio Cichella, Isaac Kaminer, Claire Walton, and Naira Hovakimyan.
\newblock Optimal motion planning for differentially flat systems using bernstein approximation.
\newblock {\em IEEE Control Systems Letters}, 2(1):181--186, 2017.

\bibitem{cichella2019consistent}
Venanzio Cichella, Isaac Kaminer, Claire Walton, Naira Hovakimyan, and Ant{\'o}nio~M Pascoal.
\newblock Consistent approximation of optimal control problems using bernstein polynomials.
\newblock In {\em 2019 IEEE 58th Conference on Decision and Control (CDC)}, pages 4292--4297. IEEE, 2019.

\bibitem{cichella2020optimal}
Venanzio Cichella, Isaac Kaminer, Claire Walton, Naira Hovakimyan, and Antonio~M Pascoal.
\newblock Optimal multivehicle motion planning using bernstein approximants.
\newblock {\em IEEE transactions on automatic control}, 66(4):1453--1467, 2020.

\bibitem{Freudenthaler2016a}
G.~Freudenthaler and T.~Meurer.
\newblock Pde-based tracking control for multi-agent deployment**financial support by the german research council (dfg) in the project me 3231/2-1 is gratefully acknowledged.
\newblock {\em IFAC-PapersOnLine}, 49(18):582--587, 2016.
\newblock 10th IFAC Symposium on Nonlinear Control Systems NOLCOS 2016.

\bibitem{FREUDENTHALER2020_a}
Gerhard Freudenthaler and Thomas Meurer.
\newblock Pde-based multi-agent formation control using flatness and backstepping: Analysis, design and robot experiments.
\newblock {\em Automatica}, 115:108897, 2020.

\bibitem{gilbert1988fast}
Elmer~G Gilbert, Daniel~W Johnson, and S~Sathiya Keerthi.
\newblock A fast procedure for computing the distance between complex objects in three-dimensional space.
\newblock {\em IEEE Journal on Robotics and Automation}, 4(2):193--203, 1988.

\bibitem{hammond2023path}
Maxwell Hammond, Venanzio Cichella, Amirreza~F Golestaneh, and Caterina Lamuta.
\newblock Path planning for continuum rods using bernstein surfaces.
\newblock {\em arXiv preprint arXiv:2312.12333}, 2023.

\bibitem{kielas2019bebot}
Calvin Kielas-Jensen and Venanzio Cichella.
\newblock Bebot: Bernstein polynomial toolkit for trajectory generation.
\newblock In {\em 2019 IEEE/RSJ International Conference on Intelligent Robots and Systems (IROS)}, pages 3288--3293. IEEE, 2019.

\bibitem{kielas2022bernstein}
Calvin Kielas-Jensen, Venanzio Cichella, Thomas Berry, Isaac Kaminer, Claire Walton, and Antonio Pascoal.
\newblock Bernstein polynomial-based method for solving optimal trajectory generation problems.
\newblock {\em Sensors}, 22(5):1869, 2022.

\bibitem{MacLennan2019_a}
Bruce~J. MacLennan.
\newblock Continuum mechanics for coordinating massive microrobot swarms: Self-assembly through artificial morphogenesis.
\newblock 2019.

\bibitem{MEURER2011_a}
Thomas Meurer and Miroslav Krstic.
\newblock Finite-time multi-agent deployment: A nonlinear pde motion planning approach.
\newblock {\em Automatica}, 47(11):2534--2542, 2011.

\bibitem{Qi2015a}
Jie Qi, Rafael Vazquez, and Miroslav Krstic.
\newblock Multi-agent deployment in 3-d via pde control.
\newblock {\em IEEE Transactions on Automatic Control}, 60(4):891--906, 2015.

\bibitem{Rubin_Book}
Miles~Barton Rubin.
\newblock {\em Cosserat Theories: Shells, Rods and Points}.
\newblock 0925-0042. Springer Dordrecht, 1 edition, 2000.

\bibitem{Selivanov2022_a}
Anton Selivanov and Emilia Fridman.
\newblock Pde-based deployment of multiagents measuring relative position to one neighbor.
\newblock {\em IEEE Control Systems Letters}, 6:1--1, 04 2022.

\bibitem{Wang2012a}
Jin-Liang Wang and Huai-Ning Wu.
\newblock Leader-following formation control of multi-agent systems under fixed and switching topologies.
\newblock {\em International Journal of Control}, 85(6):695--705, 2012.

\bibitem{Xiao2009b}
Feng Xiao, Long Wang, Jie Chen, and Yanping Gao.
\newblock Finite-time formation control for multi-agent systems.
\newblock {\em Automatica}, 45(11):2605--2611, 2009.

\bibitem{Xuan2002a}
Ping Xuan and Victor Lesser.
\newblock Multi-agent policies: From centralized ones to decentralized ones.
\newblock pages 1098--1105, 01 2002.

\end{thebibliography}


\end{document}